# A Renunciation of Axioms – A New Conception of Logic

*Eitan Wander*

*Department of Physics, Weizmann Institute of Technology, Rehovot, Israel*
*eitan.wander@weizmann.ac.il*
*ORCID 0009-0003-2496-4066*

Mathematicians invented Mathematics to escape from words, but at last they depend on them just as much as everybody else. At the end, all basic definitions will be reliant on words, yet the mathematician believes that he's elevated from them by use of axioms, only that just postpones the problem and later becomes more serious Gödelian problems due to the countability of the axiomatic array.
I think we should do our best to rid ourselves of axioms, and in this paper, I revisit what's called "Naïve Set Theory", which is set theory that is fully reliant on verbal definitions, and resolve the problems that were once found with it in the form of paradoxes, most notably, Russell's Paradox. I believe that by this, a new approach to mathematics as a whole is presented, an approach that refers to mathematics as the science of definitions .
The question "how do you define definition?" will not be a part of Mathematics but of Metamathematics and to it many of the problems known about Naïve Set Theory will be funneled.
A thorough yet concise treatment of the matter will be given in this paper.

## INTRODUCTION

We'll start with the two basic definitions in mathematics that all of the following definitions stem from.

*A point…*
is a definition that is singular.
*A set…*
is a selection of defined points.

A singular definition is one that cannot ever be applied to more than that one thing it's applied to. The definition of chair is not singular, it needs to be the definition of a specific chair in order to be singular. You can call the chair whatever you wish but it must only apply to that one single chair. For that matter you can't invent a new type of chair, say you call it a mesotropic chair, produce one instance of it, and say it's singularly defined, since you can make more mesotropic chairs in the future. The term point is used interchangeably with object.

Some variants of set theory use what's called atoms, these are objects that aren't sets (or at least treated as if they are basic) and that are differentiable from one another. Under this conception of logic the existence of these objects is allowed, and they might even provide a stern bridge between mathematical abstractness and everyday life and definitions. After all, there is no reason why a chair couldn't be used as such an atom. Nevertheless, the Von-Neumann construction of set theory using nothing more than the empty set is still perfectly allowed.

We're going to define the entirety of mathematics using these definitions and nothing else. Everything will be proven by definition and without axioms. Historically, the first conception of logic and set theory was also by definition, but it had many problems in the form of paradoxes that made researchers concede to axioms that tell you what you can and cannot do. This form of set theory is called today Naïve Set Theory [7]. We're going to revisit this approach and show how these paradoxes could be resolved, that is, why they aren't paradoxes in this conception due to a misstep along the way. A taste of this will be given right in this introduction in the resolution of the famous Russell's Paradox [7] and a treatment of other paradoxes will be given in a dedicated section.

Russell's Paradox states:
Let A be the set of all sets that don't belong to themselves.
A couldn't belong to itself since all sets in A don't belong to themselves.
But if it doesn't belong to itself then it must belong to A, that is, itself.

A set is said to be a selection of defined points, so the definition of A means that A is the set of all already defined sets that don't belong to themselves. That is, A doesn't belong to itself, not because it's not a set, but because it wasn't defined when A selected all sets. So, the moment you define a set of all sets, it's no longer a set of all sets, it can only be a set of all sets apart from himself. The set of all sets that could possibly be defined is not a set simply because a set must be a selection of defined points, and not *possibly defined points.*

This approach is called spaciotemporal definition. It means thinking of definition as dependent on time and space. The dependence on time is clear from the example; the dependence of space is because a set of all sets, for example, generally won't be the same set for two different people since they haven't defined the same sets up to that point. The definition space is the conceptual space where all singular definitions lie. It starts from being empty and every new singular definition expands it. There could be many different definition spaces.

This is a central paradox that was the first and most prominent nail in the coffin of naïve set theory, and we just saw how it's not a paradox at all in this conception. A few canonical examples are shown in the dedicated section that give a comprehensive overview of how all paradoxes of naïve set theory can be resolved. A conception which has no paradoxes is called a *consistent conception*. By that comprehensive overview, consistency is claimed.

The problem still rests, "what's going to tell us what is possible and what is impossible to do?". Once we delve into this problem, it's clear that it's not the vagueness of the definition of what a set is or what a point is that is problematic, it's the illusive "defining what a definition is" that is problematic. We'll see that the fact that you can select defined points to make any set you wish isn't problematic, and common sense will take us pretty much as far as we want to go with it. The problem is that it's hard to figure out what could be defined in the first place. That in turn is only hard because we don't really know what the meaning of definition is. The question, "What is the definition of a definition?" is clearly paradoxical and couldn't be answered in a straightforward manner, or more precisely, in a mathematical manner.

*Mathematics…*
is the theory of definitions.

We see that the question at hand couldn't even be considered as a mathematical question to begin with. An attempt of a treatment of it in mathematical terms will prove itself just as paradoxical as the initial question itself. This problem is a part of the subject of Metamathematics, in fact, it should be considered the fundamental problem of Metamathematics.

A problem in metamathematics couldn't be solved by mathematical means, so what have we to deal with it? I believe the question at hand should be dealt with in philosophical and perhaps linguistic measures, by Aristotelian Logic. I'll provide a review of the way we should go about this endeavor in a subsequent subsection. While that subsection will be metamathematical, the rest of this paper is purely mathematical.

It's important to note that the definition paradox is just as existent in the classic conception of logic as it does in this conception, yet again, simply postponed using axioms. Unsurprisingly, the axioms that claim anything about this problem are the most controversial axioms in ZFC – the axioms of powerset and choice. Both require the ability to define a certain type of definition. We'll call these axioms definition axioms and the others set axioms. We'll see that all set axioms of ZFC arise trivially from the definition of a set and a point, yet the produce of the definition axioms from the basic definitions is non-trivial.

The first chapter in this paper goes through the widely accepted ZFC axioms and provides an explanation for why they're true in this conception of logic. The second chapter visits a few canonical paradoxes that are supposedly problematic in naïve set theory. The third chapter provides a thorough introduction into the needed metamathematical bedrock. The fourth chapter and fifth chapter discuss the implications of such a conception on the Continuum Hypothesis and Gödel's Incompleteness Theorem.

## I. ZFC SET AXIOMS

The set axioms of ZFC [8] arise quite trivially from the definitions of set and point. We'll go through each one and prove it. A treatment of the definition axioms will be given in the section concerning metamathematics.

*Extensionality*
Two sets that select the same points are equal
*Proof*
Trivial, equality within a definition space is a symbol stating that two other symbols are referring to the same definition.

*Foundation*
There exists no infinite series of sets $a_n$ such that
$$... \in a_2 \in a_1 \in a_0$$
*Proof*
You could only define sets as selection of defined points, and at first there are no sets that are defined, so you simply can't construct such a descending series of sets.

*Pairing*
For any two sets, there exists a set that contains the two sets and only them as elements
*Proof*
Trivial, both sets are objects and the set of both of them is definable

*Union*
Provided a set of sets A, there exists a set containing all points selected by the elements of A
*Proof*
Trivial, select all points that are selected by the elements of A

*Infinity*
There exists an infinite set
*Proof*
There is no problem defining infinitely many points if they're singularly defined, for example, the von-Neumann ordinals up to $\omega$ are definable and provide the set of natural numbers.

*Specification*

If A is a set, $\varphi$ is a formula in the language of ZFC, and $\omega_1, \ldots, \omega_n$ are free variables, then B is a set such that
$$\forall_x\ x \in B \Leftrightarrow \big(x \in A \wedge \varphi(A, x, \omega_1, \ldots, \omega_n)\big)$$
*Proof*
Trivial, you are free to select points using formulas

*Replacement*
If A is a set, $\varphi(A, x, y, \omega_1, \ldots, \omega_n)$ a formula in the language of ZFC where x is any member of A and $\omega_1, \ldots \omega_N$ are free variables, then
$$\forall_{x \in A} \exists!_y\ \varphi \Rightarrow \exists_B \forall_{x \in A}(x \in A \Rightarrow y \in B \wedge \varphi)$$
*Proof*
Similarly to specification, this is trivial since you're free to use functions to select points for a set.

## II. A TREATMENT OF PARADOXES

*The Burali-Forti Paradox* [2]
Supposing the existence of a set of all ordinals leads to a contradiction. Since the set itself is an ordinal, it must be less than itself by the membership relation in contradiction to it being itself.
The resolution of this paradox comes together with the important definition of practical definability.

*Practical Definability…*
is the ability to define a term in a singular manner.

In other words, saying something isn't practically definable is not a statement about the philosophical meaning of being able to define such a set, it is beforehand, the ability to state what it is that you want to define in an unequivocal manner; it's not an abstract statement on the validity of such a definition, it's a practical one regarding one's ability to state the definition in the first place.
The resolution of this paradox comes from the fact that all ordinals couldn't be practically defined. If you could have defined them as points, you could have selected them to a set, but you can't define them in the first place. In a loose manner, this is because any sort of generalization you can make to construct all ordinals could be generalized further in what surmounts to be a bigger ordinal. In fact, the nature of ordinals is exactly that which can't be generalized in context of the membership relation.
The construction of the ordinals goes as follows:
1. Take the empty set {} to be 0
2. Take $n + 1$ to be $\{0, \ldots, n\}$
3. Proceed by induction…

Only then you notice that you can take the union of all the natural numbers with zero and get a set $\omega$ which is also an ordinal. Then you have
4. $\omega + n + 1 \stackrel{\text{def}}{=} \{\omega + 0, \ldots, \omega + n\}$
5. Proceed by induction…

You might think that's it but then you realize you can have $2\omega$ which is the union of all $\omega + n$. You could generalize that and get $m\omega + n$ and generalize that further to get quadratic equations and then all polynomials, but that's not the end either since you can select all polynomials to a set and get the tetration of $\omega$ by $2 - \omega^\omega$. Here you might assume that you can get away with generalizing over all hyperoperations – addition, multiplication, tetration, … - but again, you can select them all to a set and start adding ordinals all over again, as you can do indefinitely.
The Burali-Forti paradox in this new conception of logic isn't a paradox at all, it is a proof by negation for the practical indefinability of the ordinals. This is a good point to recount how paradoxes are an unnatural thing, and if your system of logic has paradoxes in it, it's wrong. By resolution of a paradox, the meaning is for an attendance to a form of doubt people might raise to this new conception of logic by means of a paradox. Hopefully, all paradoxes that might be raised could be resolved.
The infinity of ordinals, like the infinity of cardinalities, or of possible sets, is not an infinity that we're used to. It's not an infinity that we could talk about very much further than the fundamental statement that if you think you've realized them all, you haven't. This infinity will be called *Impractical Infinity*.
A very similar treatment could be given to Cantor's Paradox that states there couldn't be a set of all cardinalities. The resolution of this paradox will be that the amount of cardinalities is impractically infinite.
Such paradoxes whose resolution proves instead the impractical definability of a set will be called *Naïve Paradoxes* in retaliation for the demeaning term "naïve" set theory.

*König's Paradox* [6]
Given that the real numbers are well-ordered, and that there is a countable amount of "finitely definable" real numbers, there must be a first real number that isn't finitely definable. This provides a definition for that number in contradiction to it being non-finitely definable.
There are a few things to unpack in this paradox considering the new conception. First, what is meant by finitely definable couldn't have anything to do with definition since all the real numbers are presumed to be already defined otherwise you wouldn't be able to talk about them let alone define a subset of them. What is meant by finitely definable is what we'll call *finitely nameable*.

*Nameability…*
is the ability to assign a symbol to a definition in a non-arbitrary fashion, such that in two distinct definition spaces, constructed identically, these definitions will be *cross-equal*.

Cross-Equality is an attribute of definitions in two separate definition spaces that are the same object up to the fact that atoms defined identically are also cross-equal. If you take a real number 'Let x be a real number', then you haven't named it since that is an arbitrary definition. Another person, with a different definition space, but with the same exact construction as this, will inevitably pick a different real number to be x. More will be said about arbitrary definitions in the metamathematics section.

What is referred to in this paradox has to be finite-nameability, the ability to name a definition in a finite construction. This must be in a non-arbitrary fashion as stated in the definition of nameability, since the possibility of arbitrarily naming a definition exists for all objects, so all definitions are nameable to that extent and the set of such definitions has no reason of being countable.

The resolution of the paradox is now clear, the naming of the "first real number" is arbitrary since the construction of a well-ordering of the reals is arbitrary by use of the axiom of choice, so it isn't really named at all. The arbitrariness of the axiom of choice will also be referred to in the section about metamathematics.

*Curry's Paradox* [3]
This paradox is somewhat of a generalization of Russell's. It states that $X \stackrel{\text{def}}{=} \{x \mid x \in x \rightarrow Y\}$ is a set if and only if Y is true, where Y is a predicate. The proof goes as follows.

1. $X = \{x \mid x \in x \rightarrow Y\}$
   By definition
2. $x = X \rightarrow (x \in x \leftrightarrow X \in X)$
   Substitution
3. $x = X \rightarrow ((x \in x \rightarrow Y) \leftrightarrow (X \in X \rightarrow Y))$
   Addition of a consequent (2)
4. $X \in X \leftrightarrow (X \in X \rightarrow Y)$
   Law of concretion (1 and 3)
5. $X \in X \rightarrow (X \in X \rightarrow Y)$
   Biconditional elimination (4)
6. $X \in X \rightarrow Y$
   Contraction (from 5)
7. $(X \in X \rightarrow Y) \rightarrow X \in X$
   Biconditional elimination (4)
8. $X \in X$
   Modus ponens (6 and 7)
9. Y
   Modus ponens (8 and 6)

This line of reasoning must be false at some point since worse than the conclusion (9) it supposedly proves statement (8) which can't be true since $X = \{x \mid x \in x \rightarrow Y\}$ is simply the set of all predefined sets because $(x \in x) = False$, and we already know that $X \notin X$.

The logical fallacy of this paradox when stated in our conception might be better hidden, yet it's the same as the one in Russell's Paradox. It takes place at line (4), where it's assumed, understandably, that $(x \in x) \rightarrow Y$ from (3) produces $x \in X$, but this isn't true since x, being identically X, wasn't defined when X was, so even though $(x \in x) \rightarrow Y$ it doesn't imply that x is in X.

By these examples of common paradoxes found in naïve set theory, it is clear how all paradoxes could be resolved. Either by view of spaciotemporal definition, as with Russell's and Curry's paradoxes, or by impractical definability, as with the Burali-Forti Paradox, or by the distinction of definability and nameability.

## III. METAMATHEMATICS

The question "how do you define definition?" is paradoxical since there is no way of knowing whether you've answered it or not. A better question would be "what is the meaning of the word *definition*?" and there are a few ways of looking into it. People commonly say something is defined, or properly defined, if they can find a description that is precise, so that it can only be applied to it and nothing else. To the entirety of what you mean to say and not anything more than that. So, does that mean that if you can define something precisely than it simply could be defined? This is called the *precision hypothesis*.

In this approach there's a shift of focus from the usual disposition that asks, "what are we allowed to do?" to another disposition that asks, "what aren't we allowed to do?". A shift from a "Why?" to a "Why Not?".

The line of reasoning that begins with "People commonly say…" is an Aristotelian form of logic that enquires into the true meaning of words. *Vox Populi, Vox deis* – the voice of the people is the voice of God. In his book "Nichomachean Ethics" [1] Aristotle has used this form of logic in order to enquire into the meaning of abstract terms such as "virtue," "happiness," "justice," and "goodness," among others. This is the method I think should be implemented to the enquiry of the meaning of definition.

Naturally, the voice of the majority isn't always right, but the method of *vox populi* doesn't call for a democratic election over each statement for it to be true. On the contrary, some statements might be considered true even though a majority would have elected them false.

The method of *vox populi* starts with as little and as least controversial statements that are to be suggested as true since they are the voice of the people, and then shall be used axiomatically in a strict logical manner. So a person might elect a false statement if it were shown to him out of the blue, but in theory the method of *vox populi* states that if that person was logical and he were to follow and understand the steps of the proof perfectly, at the end of reading it he would vote the other way around.

When you come to think about it, we've been using Aristotelian logic all along in this paper. In fact, the method of *vox populi* is the natural method of interrogating the meanings of abstract terms.

This isn't to say that any argument a person gives with a veil of Aristotelian logic – for instance, the writer of this paper – must be true. The simple fact that someone says "people commonly say…" and then uses dubious logic that people commonly agree with to get where he wants doesn't make it true, and there could be many contradictions given to a person that abides by such an ideology. Aristotelian logic is really about conversation, and it always remains somewhat open-ended. At the end of the day, if a considerable percentage of logical individuals follow your reasoning and still don't agree with you, then your argument is wrong. The following lines of reasoning I apply in this chapter aren't entitled (by the author at least) to be of a finality of any sort. While the treatment of the basis of Mathematics in this paper is proposed to be final, the treatment of Metamathematics isn't. The following lines are a proposed resolution by the author and are mainly suggested as a spark of conversation.

Aristotelian logic suggests the precision hypothesis is true, which means a definition is allowed if it is *precise* and *exact*, yet that doesn't mean that all definitions are singular definitions. The precision hypothesis allows the existence of *descriptive definitions*, definitions that describe other definition, i.e. every set with a binary function that obeys the three axioms is called a *group*, but the definition of a group certainly isn't singular.

Does the precision hypothesis predict the correctness of the powerset axiom? For the powerset axiom to be true we need to be able to define all subsets of a given set at once and then select them to a set. The question now rests, is a subset of a set a singular definition or a descriptive definition? It could be thought of as descriptive, since you say a subset of a set is a set that selects objects that the parent set also selects. But then again, here you're talking about the subsets of a singular set. Let the set be A, then 'a subset of A' is also a descriptive definition, so seemingly the powerset axiom couldn't be correct. The only thing is, you're perfectly capable of uttering sentences the likes of "define B to be *a* subset of A", now what kind of definition is that for the set B? Is it descriptive or singular? Well, it must be singular, since there isn't any other thing other than the set B that you could claim to be B, but this still seems a little fuzzy, since you haven't really *chosen* B. This is what's called an *arbitrary definition*, and an arbitrary definition is indeed singular.

Again we are confronted with a question that we instinctively attempt to barge at with a "why?", when the truth of the matter is "why not?". Why couldn't we define a set arbitrarily? Why wouldn't it be singular if we can check whether its' items are identical to any other set even though they were chosen arbitrarily?

But maybe you can't *really* define arbitrary definitions, maybe they shouldn't be allowed. How are we supposed to find what kinds of definitions are we allowed to do and what kinds we aren't? Well, using the precision hypothesis, we should enquire as to whether the definition of B in this case is precise.

What do we *mean* by B? Is there any other thing that can be confused with B? What are the properties of B?

Since any other subset of A can be quickly reassured to not be B, and any other object or definition cannot be B since it isn't a subset of A, we have that the definition is precise and therefore there is no reason why we shouldn't be able to define it.

So, you can define any subset of any set and therefore you can also define all of them, because each one is again singularly defined, and you can select them all to a set so that you have the powerset. And thus, the powerset axiom is obtained.

Now what about the axiom of choice? The treatment is similar, you can define an arbitrary set that selects one object that each set in a superset selects since it's a precise definition.

Should arbitrary definitions be possible? Say you have a set, and you *define* all subsets of that set, it *feels* a bit as if you don't really know all the sets you have defined. Do you have to know everything that you have defined? Well there is really no reason not to, take for instance the natural numbers, you don't really know every number personally, you haven't gone through them one by one and defined them yourself, you used recursion and the infinity axiom. It comes down to your ability to assert whether something was defined or not defined. Say you define all subsets of a set, and then you are presented with one, it's perfectly reasonable that you'll be able to look at it and say, "this is a subset of set A, which I have defined". Perhaps that should be the ultimate test of definition, whether when presented with a singular object you can say, "yes, I have defined it".

What would be an Aristotelian logical backing to such a statement? Well, it is common to say that something has a certain property if you can assert that it has it. Meaning, if you saw a box that was green, it's only logical that it is a green box. So it seems to be an accurate test. It's also reasonable to say that if you look at something and say, "it is singular, and I have defined it" then it simply is an object.

A contrary example to this is if you try and define all ordinals, you could perfectly well test if a certain set is an ordinal once presented with one, but you can't define all ordinals as was shown above. So if you can say "define all ordinals" and then look at every set and test whether it's an ordinal or not, then why can't you define all ordinals? The resolution of this comes from the fact that you could be presented with ordinals that select objects that you haven't defined, and so you can't possibly say that said ordinal is defined since it's not even a set. For instance, you are permitted to define all ordinals that could be defined as sets of predefined objects.

Another noteworthy school might suggest that a definition could be defined simply if it doesn't resolve in a paradox. In other words, if you can't prove it not to be a definition using a paradox, then it is a definition. This school will be called *the paradox school* in contrast to the *Aristotilean School*.

## IV. THE CONTINUUM HYPOTHESIS

The continuum hypothesis [4] states that there is no cardinality between that of the natural numbers and that of the reals. It was demonstrated to be independent of the ZFC axioms, what does that mean in this new conception?
The fact that CH (Continuum Hypothesis) Is Independent of ZFC means that you can't construct a set with a cardinality between the naturals and the reals using the powerset and choice definition axioms. In our conception this means that CH is true if and only if there isn't a definition that couldn't be constructed using powerset and choice.
This reduces CH to a metamathematical problem of construction. Is there a set that you could define that isn't constructed using the powerset axiom or the axiom of choice?

## V. GÖDEL'S INCOMPLETENESS THEOREM

Gödel's theorem [5] roughly states that in any consistent axiomatic system there are statements that can neither be proved nor disproved. This theorem doesn't apply to this conception since it isn't axiomatic, question is, are there similar problems in this conception?
A statement is true in mathematics under this conception if and only if it is true by definition. Something is true by definition if it's either plainly obvious or can be deduced logically from plainly obvious statements. Such plainly obvious statements (POS) could be thought of as an axiomatic system to which Gödel's theorem could be applied. It doesn't apply, however, if the amount of POS is infinite in a way that is non-recursively enumerable.

That is, couldn't be laid down by a finitely written algorithm.
Are our basic definitions' POS recursively enumerable? Say we had a finitely written algorithm, written verbally, that produces all the POS for the point and the set.
This is equivalent to having an algorithm that determines whether a statement is plainly obvious or not plainly obvious. This is because we can enumerate any sentence in English using ASCII coding and create an algorithm that runs through all possible strings and pass them through such an algorithm.
This way, Gödel's theorem is reduced to a fundamental question in Aristotelian logic, "Is human affirmation of statements as plainly obvious or otherwise computable?"

In the paradox school the whole concept of POS is irrelevant, so the incompleteness theorem and its' equivalencies remain untrue for this school in the meanwhile.

## VI. CONCLUSION

With this conception of logic, a renunciation of axioms is given. This is an exciting thing since it has the potential of bridging the gap between abstract mathematics and everyday life.

In the tree of knowledge, where the trunk is the base suppositions and the branches are the theorems presented thereupon, Mathematics claimed to be the form of knowledge built from the ground up. Now we see that this tree has roots, and while the trunk is singular, the roots very much aren't. All of Mathematics was built on words, here all that was given is a thorough treatment of this already present yet unnoticed or repressed truth.
We see that while the ground is firm, and we still have all the solid earth we need to continue exactly as we have been, expanding mathematics rigorously, we also need to consider the growth of the tree of knowledge *downwards* – the growth of its' roots into the soil.
The reason this growth is inside the earth and could only increase the strength and balance of the tree is because mathematics is *true*. There isn't anything we can say that will shake its foundations when we already know it to be useful. Since Aristotelian logic is the logic of common sense no one will ever renounce what is useful since to a large extent there is no greater truth than that which is useful and sturdy enough to hold on to.
All inquiry into Metamathematics, the realm beneath the surface of the ground, could only strengthen our assertion that mathematics is *true* with more and more compelling and grounded argumets to why we should be able to do what we are doing with definitions.

Mathematics is the theory of definitions, and we come upon definitions every time we utter a word. I believe that with a ground enough understanding of Metamathematics a new and exciting branch of mathematics could be realized, the Mathematics of *common words*. From such an achievement we might finally be able to bridge the gap and obtain absolute true and false into the arguments of everyday life. An achievement that's now only feasible in "lab conditions" with the two basic definitions of point and set.

A system which asserts the truth of a statement presented in everyday words would be revolutionary for the development of humanity. These days, we are all much to aware of the effects of populistic arguments and exposure through social media to content that only serves to cement the individual's sitting in his stance and further polarize our society. Humanity has thought, naively, that if people have the information, they will reach the right conclusions, yet without taking any sides, two contradicting conclusions can't possibly be right at the same time. I believe in a more structured presentation of reality. One which gives statements as true even though people might not agree with them. One whose truth is given from the fact that under the common (albeit logical) individual's further inquiry into the argument that ended with this statement an acceptance of it as true will be achieved by the individual.

## REFERENCES

[1] Aristotle, Nicomachean ethics. Hackett Publishing, 2019.

[2] Burali-Forti, Cesare. "Una questione sui numeri transfiniti." Rendiconti del Circolo Matematico di Palermo (1884-1940) 11.1 (1897): 154-164.

[3] Curry, Haskell Brooks. Foundations of mathematical logic. Courier Corporation, 1977.

[4] Gödel, Kurt. "Consistency of the continuum hypothesis." (1940).

[5] Gödel, Kurt. "On formally undecidable propositions of principia mathematica and related systems i 1 (1931)." Godel's Theorem in Focus. Routledge, 2012. 17-47.

[6] König, Julius. "Über die Grundlagen der Mengenlehre und das Kontinuumproblem." Mathematische Annalen 61.1 (1905): 156-160.

[7] Russell, Bertrand. "Mathematical logic as based on the theory of types." American journal of mathematics 30.3 (1908): 222-262.

[8] Zermelo, Ernst. "Untersuchungen über die Grndlagen der Mengenlehre. I." Mathematische Annalen 65.2 (1908): 261-281.